\documentclass{amsart}
\usepackage{amsfonts,amscd,amssymb,amsmath,amsthm}
\usepackage{graphicx,mathrsfs,appendix}
\usepackage{centernot,caption,subcaption,color,url}
\usepackage[T1]{fontenc}
\usepackage{ctable}

\begin{document}
\newtheorem{theorem}{Theorem}
\newtheorem{lemma}[theorem]{Lemma}
\newtheorem{conjecture}[theorem]{Conjecture}
\newtheorem{cor}[theorem]{Corollary}
\newtheorem{proposition}[theorem]{Proposition}
\newtheorem{claim}[theorem]{Claim}
\newtheorem{remark}[theorem]{Remark}
\newcommand{\R}{\mathbb{R}}
\newcommand{\W}{\mathcal{W}}
\newcommand{\T}{\mathcal{T}}
\newcommand{\B}{\mathcal{B}}
\newcommand{\A}{\mathcal{A}}
\newcommand{\G}{\mathcal{G}}
\newcommand{\F}{\mathcal{F}}
\newcommand{\Z}{\mathbb{Z}}
\newcommand{\Q}{\mathbb{Q}}
\newcommand{\E}{\mathbb E}
\newcommand{\N}{\mathbb N}
\newcommand{\1}{\mathbf 1}

\title[Spatial preferential attachment with choice-based edge step]{Spatial preferential attachment with choice-based edge step}
\author{Yury Malyshkin}
\address{Tver State University}
\email{yury.malyshkin@mail.ru}
\subjclass[2010]{05C80}
\keywords{random graphs, preferential attachment, power of choice, fitness}
\date{\today}

\begin{abstract}
We study the asymptotic behavior of the maximal in-degree in the spatial preferential attachment model with a choice-based edge step. We prove different types of behavior of maximal in-degree based on the model's parameters.
\end{abstract}


\maketitle

\section{Introduction}

In the present work, we study the addition of the choice-based edge step to the spatial preferential attachment model. Preferential attachment models are widely used to describe complex networks (see, e.g., \cite{Hof16}). The general idea behind preferential attachment models is that vertices with higher degrees are more likely to attract edges from newly introduced vertices.
In the classical preferential attachment graph model (\cite{BA99}), vertices are undistinguished and the probability of drawing an edge to a vertex is proportional to its degree's linear function (see, e.g., \cite{Mori02, Mori05}). The spatial preferential attachment model adds geometry to this procedure by giving each vertex a coordinate (see, e.g. \cite{ABCJP09}). In this model, each vertex has a neighborhood of a size proportional to a linear function of its degree, and we draw an edge to it if a new vertex belongs to such a neighborhood.

We modify this model by adding a choice-based edge step to the model to allow connection between far away vertices. The edge step is used to draw the edges between the old vertices (see, e.g., \cite{ARS21}), which in our case could be far away vertices. The choice step uses an auxiliary sample of vertices chosen by preferential attachment rules, from which we chose the vertex based on additional rule (see, e.g., \cite{HJ16,KR14,Mal18,MP14}). Unlike previous models, we would use two auxiliary samples and pick multiple vertices from the final sample at the same time. The addition of choice often results in the effect of condensation, when a single vertex has a linear (over the total number of edges) degree (see, e.g., \cite{HJY20,Mal20,MP15}). Picking multiple vertices at the same time allows us to get that effect on multiple vertices. 

Let us introduce our model. Fix $k,\in\N$, $d,M\in\N,$ $a,\alpha\in(0,1/2)$, $b,\beta>0,$ and distribution of a random variable $m$ with values in $1,...,M$, which are parameters of our model. We consider a sequence of graphs $G_n$. We consider i.i.d. random variables $X_1,X_2,...$ distributed uniformly in $[0,1]$, such that $X_i$ corresponds to $v_i$ and represents the location of a vertex. We also consider i.i.d. random variables $m_i$, $i\in N$, which are copies of $m$. We start with the initial graph $G_0$, which consists of $n_0>M$ vertices.  To build a graph $G_{n+1}$ from $G_n$, we add a vertex $v_{n+1}$ with location $X_{n+1}$ and draw edges in two steps.
\begin{itemize}
\item[Vertex step:] we draw edges from $v_{n+1}$ to the vertices that satisfy $X_{n+1}\in B_n(v_i)$, where
\begin{equation}
    B_n(v_i)=\left\{x:|x-X_i|<\frac{a/2\deg^{-}_{G_n}v_i+b/2}{n}\right\},
\end{equation}
where $|\cdot|$ is the torus norm in $[0,1]$ and $\deg^{-}_{G_n}v_i$ is the in-degree of $v_i$ in $G_n$.
\item[Edge step:] We choose the vertex $u_n$ uniformly among all vertices. We consider $d$ samples $y_n^{i,1},y_n^{i,2},...$, $i=1,...,d$  of vertices of $G_{n}\cup v_{n+1}-u_n$, chosen independently with probabilities $\frac{\alpha\deg^{-}_{G_n} v+\beta}{n}$. From each sample, choose $m_n$ vertices with the highest in-degree in the sample to a secondary sample $z_n^{1},...,z_n^{dm_n}$ (the sample would contain enough vertices for large enough $n$; otherwise, choose the rest uniformly among all vertices). Then, we draw edges from $u_n$ to $m_n$ different vertices from the sample with the highest in-degree (in the case of a tie, chosen randomly, it would not affect the degree distribution).
\end{itemize} 

\section{Results}
Let us formulate our main results. Let $M_i(n)$ be the $i$-th highest in-degree of the vertices of $G_n$. Consider the evolution of $M_i(n)$. On the $n+1$-th step, it could be increased in two ways.

First, $v_{n+1}$ could be in the neighborhood of the vertex with in-degree $M_i(n)$. The probability (conditioned on graph $G_n$) to do so is at least (exactly if there is a single vertex with the highest in-degree) $\frac{aM_i(n)+b}{n}$.

Second, we could draw an edge to it as the second vertex during an edge step. To do so, we need a vertex of in-degree $M_i(n)$ to appear in the sample and have at most $m_n-1$-th highest in-degree vertices in the sample. For any vertex $u$, the probability of appearing in the initial sample is
$$\frac{\alpha\deg^{-}_{G_n} u+\beta}{n}$$
and the probability of appearing in at least one of the initial samples equals
$$1-\left(1-\frac{\alpha \deg^{-}_{G_n} u+\beta}{n}\right)^{d}.$$
Let us introduce a function
$$h(x):=1-\left(1-\alpha x\right)^{d}$$
Let $p_i(n)$ be the probability for the $i$-th highest in-degree vertex to appear in one of the initial samples. We get 
$$p_i(n)=h\left(\frac{M_i(n)}{n}+\frac{\beta}{\alpha n}\right)=h\left(\frac{M_i(n)}{n}\right)+O\left(\frac{1}{n}\right).$$  
If the vertex has the highest in-degree, it would always be chosen.
For $i>1$, to increase $M_i(n)$, we would need at least $i-m_n$ of the $(i-1)$-th highest in-degree vertices to not be in the sample and a vertex with in-degree $M_i(n)$ to be in the sample. It happens with probability (if $M_{i-1}(n)>M_{i}(n)>M_{i+1}(n)$) 
\begin{equation}
\label{eq:prob}
\begin{gathered}
p_{i}(n)=h\left(\frac{M_i(n)+\frac{\beta}{\alpha n}}{n}\right)Q_{i,m_n}\left(\frac{M_1(n)+\frac{\beta}{\alpha n}}{n},...,\frac{M_{i-1}(n)+\frac{\beta}{\alpha n}}{n}\right)\\
=h\left(\frac{M_i(n)}{n}\right)Q_{i,m_n}\left(\frac{M_1(n)}{n},...,\frac{M_{i-1}(n)}{n}\right)+O\left(\frac{1}{n}\right),
\end{gathered}
\end{equation}
for some polynomial function $Q_{i,m_n}$, whose coefficients are dependent on $m_n$ (we could consider set of non-random polynomials $Q_{i,r}$ and choose one of them dependent on the value of $m_n$, we also put $Q_{1,r}=1$). Note that 
\begin{equation}
\label{eq:polinom}
\begin{gathered}
Q_{i,r}(x_1,...,x_{i-1})= Q_{i+1,r}(x_1,...,x_{i}),\quad i\leq r-1\\
Q_{i,r}(x_1,...,x_{i-1})> Q_{i+1,r}(x_1,...,x_{i}),\quad i\geq r
\end{gathered}
\end{equation}
for any $2\geq x_1\geq x_2,...,\geq x_{i}>0$.
Let us define polynomials $Q_i(x_1,...,x_{i-1}):=\E Q_{i,m_n}(x_1,...,x_{i-1})$ and let $r_m:=\min \{r:\Pr(m=r)>0\}$. Due to \eqref{eq:polinom} we get
\begin{equation}
\label{eq:polinom_est}
\begin{gathered}
Q_{i}(x_1,...,x_{i-1})= Q_{i+1}(x_1,...,x_{i}),\quad i\leq r_m-1\\
Q_{i}(x_1,...,x_{i-1})> Q_{i+1}(x_1,...,x_{i}),\quad i\geq r_m
\end{gathered}
\end{equation}
for any $2\geq x_1\geq x_2,...,\geq x_{i}>0$.

Consider functions
$$g_i(x_1,...,x_i)=ax_i+h(x_i)Q_{i}(x_1,...,x_{i-1})\quad i\in\N.$$
Introduce functions 
$$f_i(x_1,...,x_i)=g_i(x_1,...,x_i)-x_i$$
and vector-function 
$$F_i(x_1,...,x_i)=(f_1(x_1),...,f_i(x_1,...,x_i)).$$ 
Let $K$ be the highest number that the system $F_K=0$ has a positive solution $x^{\ast}=(x_{1}^{\ast},...,x_{K}^{\ast})$. Note that $K>1$ if $a+d\alpha>1$ (since $f_1(x_1)$ is concave, $f_1(0)=0$ and $f^{\prime}_{1}(0)=a+d\alpha-1$).
Let us formulate our main theorem.

\begin{theorem}
\label{th:main}
In the model described above
\begin{enumerate}
\item if $a+d\alpha<1$, then for any $\epsilon>0$ 
$$\Pr\left(\forall n>N: n^{a+d\alpha-\epsilon}<M_1(n)<n^{a+d\alpha+\epsilon}\right)\to 1$$
as $N\to\infty$.
\item if $a+d\alpha=1$, then almost surely
$$\lim_{n\to\infty}\frac{M_1(n)\ln n}{n}=\frac{2}{(d\alpha)^2}.$$
\item if $a+d\alpha>1$, then for all $k\leq K$ almost surely
$$\lim_{n\to\infty}\frac{M_k(n)}{n}=x_k^{\ast}.$$
\end{enumerate}
\end{theorem}

\section{Maximal degrees}

In the current section, we consider maximal in-degree. The maximal in-degree always increases at the edge step if the vertex with the maximal in-degree is in the sample. Hence, if $M_1(n)>M_2(n)$, then
\begin{equation}
\label{eq:max_degree}
\E(M_1(n+1)-M_1(n)|\F_n)=g_1(M_1(n)/n)+O\left(\frac{1}{n}\right).
\end{equation}
Note that if $M_{1}(n)=M_2(n)>M_3(n)$, then probability to increase $M_1(n)$ would increase (since only one of the vertex should be preset in the sample), and probability to increase $M_2(n)$ would decrease (both vertices must be present in the sample). Hence, if we neglect this event, we would get lower bounds for $M_1(n)$ and upper bounds for $M_2(n)$.

If the process $X(n)$ with bounded increments satisfy equation \eqref{eq:max_degree}, i.e. 
$$\E(X(n+1)-X(n)|\F_n)=g_1(X(n)/n)+O\left(\frac{1}{n}\right),$$
then its behavior could be described by Lemma 3.1 of \cite{M25}. By that lemma $X(n)/n$ converges almost surely to the root of equation $g_1(x)-x=f_1(x)=0$ if $g_1^{\prime}(0)=a+d\alpha>1$, $X(n)\ln n/n$ converges almost surely to $\frac{-2}{g^{\prime\prime}(0)}=\frac{2}{(d\alpha)^2}$ if $a+d\alpha=1$ and $\Pr\left(\forall n>n_0: n^{a+d\alpha-\epsilon}<X(n)<n^{a+d\alpha+\epsilon}\right)\to 1$ as $n_0\to\infty$ if $a+d\alpha<1$. 

That would give us the required results for $M_1(n)$ if event $M_{1}(n)=M_2(n)$ does not happen after some moment. Let us consider two cases: $r_m>0$ and $r_m=0$ (meaning $\Pr(m=1)>0$). 

In the case $\Pr(m=1)>0$ due to \eqref{eq:polinom_est} the derivative at $0$ for the representation of the increment of $M_{2}(n)$ would be lower than for $M_1(n)$, and, therefore, upper bound for $M_2(n)$ would be lower than lower bound for $M_{1}(n)$, which means $M_1(n)>M_2(n)$ after some moment. 

In the case $r_m>0$, under condition $M_1()>M_2(n)>...>M_{r_m}(n)>M_{r_m+1}(n)$ increments of $M_1(n),M_2(n),...,M_{r_m}(n)$ would be described by the same equation, while function in $M_{r_m+1}(n)$ increment would have lower derivative at $0$, and, therefore, $M_{r_m+1}(n)$ would be lower after some moment. The equation \ref{eq:max_degree} would then provide lower bound for $M_1(n)$ and upper bound for $M_{r_m}(n)$, which means if $M_1(n)=M_{r_m}(n)$ infinitely many times, then they converges (after normalizing in the sense of theorem~\ref{th:main}) to the same limit, and so does all $M_i(n)$ between them. If $M_1(n)=M_{r_m}(n)$ only finitely many times, then $M_1(n)>M_{r_m}(n)$ after some moment, we could repeat the argument to $M_1(n),M_2(n),...,M_{r_m-1}(n)$ and $M_2(n),...,M_{r_m-1}(n)$ to obtain that either they intersect infinitely many times or isolate after some moment, in both cases they converge to the same limit.

\section{Case $a+d\alpha>1$}

In the current section, we prove convergence of $M_{k}(n)/n$ for $k\leq K$ for the case $2a+d\alpha>1$. We rely on stochastic approximation (see, e.g., Section 2.4 in \cite{Pem07}) to do so. We would use induction over $k$.

Let us consider the evolution of $M_i(n)$.
Note that if $M_{i-1}(n)>M_{i}(n)$, then
$$\E(M_i(n+1)-M_i(n)|\F_n)=g_i(M_1(n)/n,...,M_i(n)/n)+O\left(\frac{1}{n}\right).$$
Define process $Z_i(n)=M_i(n)/n$. Then we have representation
$$\E(Z_i(n+1)-Z_i(n)|\F_n)=\frac{1}{n+1}\left(\E(M_i(n+1)-M_i(n)|\F_n)-Z_i(n)\right).$$
Let us consider vector $W_k(n)=(Z_1(n),...,Z_K(n))$. We get (if $Z_1(n)>Z_2(n)>...>Z_k(n)>Z_{k+1}$, $k\leq K$)
$$\E(W_k(n+1)-W_k(n)|\F_n)=\frac{1}{n+1}\left(F_k\left(Z_1(n),...Z_K(n)\right)+O\left(\frac{1}{n}\right)\right).$$
Components $f_i(Z_1(n),...,Z_i(n))$ of $F_i(Z_1,...,Z_i)$ are concave functions of $Z_i(n)$ and depends on $Z_1(n),Z_2(n),...,Z_{i-1}(n)$ through $Q_{i}$:
$$f_i(Z_1(n),...,Z_i(n))=(a-1)Z_i(n)+h(Z_i(n))Q_{i}(Z_1(n),...,Z_{i-1}(n)).$$

Let us assume $Z_i(n)\to x_n^{\ast}$ almost surely for $i<k$ and prove such a convergence for $k$.
The case $k\leq r_m$ was considered in the previous section. Let's assume $k>r_m$.
Recall that $(x_1^{\ast},...,x_{k}^{\ast})$ is the positive solution of quation $F_k(x_1,...,x_k)=0$ and the other solution of that equation is $(0,...,0)$. Since $f_k(x_1,...,x_k)$ is concave as function of $x_k$, and $Q_i$ is continuous function, for any $\epsilon>0$ there is $\delta-$neighborhood of $(x_1^{\ast},...,x_{k-1}^{\ast})$, in which equation $f_k=0$ has a positive solution that belongs to $\epsilon-$neighborhood of $(x_1^{\ast},...,x_{k}^{\ast})$. For large enough $N$, all $Z_i(n)$, $i<k$, belongs to that $\delta-$neighborhood for $n>N$. Due to the definition of $Q_k(x_{1},...,x_{k})$, it decreases with increase of $(x_{1},...,x_k)$. Hence, increments of $Z_{k}(n)$ for $n>N$ would be bounded from below by
$$\frac{1}{n+1}\left((a-1)Z_k(n)+h(Z_k(n))Q_{i}(x_1^{\ast}+\delta,...,x_{k-1}^{\ast}+\delta) +O\left(\frac{1}{n}\right)\right)$$
and from above by
$$\frac{1}{n+1}\left((a-1)Z_k(n)+h(Z_k(n))Q_{i}(x_1^{\ast}-\delta,...,x_{k-1}^{\ast}-\delta) +O\left(\frac{1}{n}\right)\right).$$
Therefore, we could estimate $Z_k(n)$ by processes $X_{k}(n)=X_{k}(n,\epsilon)$ (from below) and $Y_k(n)=Y_{k}(n,\epsilon)$ (from above) with bounded increments, such that
$$\E(X_k(n+1)-X_k(n)|\F_n)=\frac{1}{n+1}\left((a-1)Z_k(n)+h(Z_k(n))Q_{i}(x_1^{\ast}+\delta,...,x_{k-1}^{\ast}+\delta) +O\left(\frac{1}{n}\right)\right),$$
$$\E(Y_k(n+1)-Y_k(n)|\F_n)=\frac{1}{n+1}\left((a-1)Z_k(n)+h(Z_k(n))Q_{i}(x_1^{\ast}-\delta,...,x_{k-1}^{\ast}-\delta) +O\left(\frac{1}{n}\right)\right).$$
Hence, due to stochastic approximation results (see, e.g., Corollary 2.7 in \cite{Pem07}), $X_k(n,\epsilon)$ and $Y_{k}(n,\epsilon)$ converges almost surely to a point in $(x_{k}^{\ast}-\epsilon,x_{k}^{\ast}+\epsilon)$. As a result, by tending $\epsilon$ to $0$, we get that $Z_{k}(n)\to x_{k}^{\ast}$ almost surely.

\section*{Acknowledgements.}
The presented work was funded by a grant from the Russian Science Foundation (project No. 24-21-00247).

\end{document}